\documentclass[12pt]{article}
\usepackage{euscript}
\usepackage{amsfonts}
\usepackage{amsthm}
\usepackage{amsmath}
\usepackage{amsfonts}
\usepackage{latexsym}
\usepackage{amssymb}

\begin{document}

\newtheorem{lemma}{Lemma}[section]
\newtheorem{theo}[lemma]{Theorem}
\newtheorem{coro}[lemma]{Corollary}
\newtheorem{rema}[lemma]{Remark}
\newtheorem{fed}[lemma]{Definition}
\newtheorem{propos}[lemma]{Proposition}

\def\bdem{\begin{proof}}
\def\edem{\end{proof}}
\def\bequ{\begin{equation}}
\def\eequ{\end{equation}}

\newcommand{\laba}{\label}

\renewcommand{\thesection}{\arabic{section}}
\renewcommand{\theequation}{\thesection.\arabic{equation}}
\newcommand{\equnew}{\setcounter{equation}{0}}

\newcommand{\la}{\langle}
\newcommand{\ra}{\rangle}
\newcommand{\berg}{A^2}
\newcommand{\deltab}{\tilde{\Delta}}
\newcommand{\hip}{{\cal A}}
\newcommand{\papa}{ H^{\infty} }
\newcommand{\disc}{ \mathbb{D} }
\newcommand{\noi}{\noindent}
\newcommand{\ov}{ \overline }
\newcommand{\tosot}{\stackrel{\mbox{\scriptsize{sot}}}{\to}}
\newcommand{\towot}{\stackrel{\mbox{\scriptsize{wot}}}{\to}}
\newcommand{\rr}{\mbox{$\rightarrow$}}
\newcommand{\om}{\omega}
\newcommand{\Om}{\Omega}

\newcommand{\dist}{\mbox{dist}\,}
\newcommand{\comb}[2]{#1 \choose #2}
\newcommand{\supp}{ \mbox{supp}\,  }
\newcommand{\diag}{\mbox{diag}\,}

\renewcommand{\L}{\mathfrak{L}}
\newcommand{\D}{\mathfrak{D}}

\hyphenation{imme-dia-te-ly}  \hyphenation{no-thing}  \hyphenation{di-ffe-rent}
\hyphenation{con-ti-nui-ty}   \hyphenation{se-ve-ral}   \hyphenation{e-qui-va-lent}
\hyphenation{pre-hy-per-bo-lic}   \hyphenation{hy-per-bo-lic}   \hyphenation{e-qui-va-len-ces}
\hyphenation{fo-llo-wing}   \hyphenation{a-na-lo-gous-ly}   \hyphenation{lo-ca-li-za-tions}


\title{A generalization of Toeplitz operators on the Bergman space}
\author{Daniel Su\'{a}rez}

\date{\mbox{}}

\maketitle

\begin{quotation}
\noindent
\mbox{ } \hfill      {\sc Abstract}    \hfill \mbox{ } \\
 \hfill \mbox{ }\\
{\small \noindent
If $\mu$ is a finite measure on the unit disc and $k\ge 0$ is an integer, we study a generalization
derived from Englis's work, $T_\mu^{(k)}$, of the  traditional Toeplitz operators on the Bergman space $\berg$,
which are the case $k=0$. Among other things, we prove that when $\mu\ge 0$, these operators are bounded
if and only if $\mu$ is a Carleson measure, and we obtain some estimates for their norms.
}
\end{quotation}

\maketitle

\setcounter{section}{0}

\section{Introduction and preliminaries}
\equnew

Let $\berg$ be the Bergman space of holomorphic function on the disc $\disc$ with respect to the normalized area measure
$dA$, and $\L(\berg)$ be the Banach space of bounded operators on $\berg$.
If for $z\in\disc$, $\varphi_z\in \mbox{Aut}(\disc)$ denotes the involution that interchanges $0$ and $z$,
the change of variables operator $U_zf = (f\circ\varphi_z) \varphi'_z$ is unitary and self-adjoint.
Here, $\varphi'_z = - K_z/ \| K_z \|$, where $K_z$ is the reproducing kernel for $z$, and  $\| K_z \|= (1-|z|^2)^{-1}$.

For $f, g, h \in \berg$, define the rank-one operator $(f\otimes g) h := \la h, g\ra f$.
In particular, if $e_k = \sqrt{k+1}\,w^k$  ($k\ge 0$) is the standard base of $\berg$, the operator
$E_k:= e_k\otimes e_k$ is the orthogonal projection onto the subspace generated by $e_k$.
Hence, for every $z\in\disc$ and $f, g\in \berg$ we have
$$\la U_zE_0U_z f, g\ra = (1-|z|^2)^{2}f(z) \ov{g(z)}.$$

So, if $d\tilde{A}(z) = (1-|z|^2)^{-2}dA(z)$ denotes the  invariant area measure on $\disc$ and
$a\in L^\infty$, the traditional Toeplitz operator $T_a$ can be written as
$$
T_a = \int_D U_zE_0U_z \, a(z)d\tilde{A}(z),$$
where the integral converges in the weak operator topology. This led Engli\v{s}  in  \cite{eng4} to consider operators
defined as above, where $E_0$ is replaced by more general operators $R$ that are diagonal with respect to the standard base
(a radial operator). Among other results, he proved that if $R$ is a radial operator in the trace class and $a\in L^\infty$,
then
$$  R_a := \int_D U_zRU_z \, a(z)d\tilde{A}(z)\in \L(\berg)  \ \mbox{ and }\   \|R_a\| \le \|R\|_{tr} \, \|a\|_\infty .
$$
Since such operator $R$ is a $\ell^1$-linear combination of the projections $E_j$, with the trace norm of $R$ given by
the correspondent $\ell_1$-norm of its eigenvalues, the above result is equivalent to
$$
T^{(j)}_a := \int_D U_zE_jU_z \, a(z)d\tilde{A}(z)\in \L(\berg)  \ \mbox{ and }\   \|T^{(j)}_a\| \le  \|a\|_\infty $$
for every integer $j\ge 0$.
We study this type of operators and a generalization $T_\mu^{(j)}$, where  $ad\tilde{A}$ is replaced by the expression
$(1-|z|^2)^{-2}d\mu(z)$, for $\mu$ a measure whose variation $|d\mu|$ is a Carleson measure. As in the well known case $j=0$,
these operators turned out to be bounded, and when $\mu$ is positive we find lower and upper bounds for their norms.
We also characterize compactness and show that these operators are norm limits of traditional Toeplitz operators.

Useful tools for our study will be the $n$-Berezin transform and the invariant Laplacian.
If $n\ge 0$ is an integer, the n-Berezin transform of $Q\in\L(\berg)$ is
\begin{eqnarray*}
B_n (Q)(z)  & \! :=\!  &   (n+1) \sum_{j=0}^n  {n \choose j} \frac{(-1)^j}{(j+1)} \la Q U_z e_j, U_z e_j\ra    .
\end{eqnarray*}
In particular, if $Q= T_\mu$, where $\mu$ is a finite measure on $\disc$, a straightforward calculation shows that
\bequ\laba{dmuda}
B_n ( \mu)(z) := B_n (T_\mu)(z) = \int_D  (n+1) \frac{ (1-|\varphi_z (\zeta )|^2 )^{n+2} }{(1-|\zeta |^2 )^2} \, d\mu(\zeta ).
\eequ
Observe that the last expression defines $B_n(\mu)$ for any measure $\mu$ of finite total variation,
even if $T_\mu$ is not bounded. In particular, if $\mu = adA$ with $a\in L^1$, we write $B_n(a):= B_n(adA)$,
which is also $B_n(T_a)$ if $T_a$ is bounded.  It is clear from the definition that $\| B_n(Q)\|_\infty  \le (n+1) 2^n \|Q\|$.
Also, it was showed in \cite{sua1} that
\bequ\laba{bobn}
B_n B_0 (Q) = B_0 B_n (Q)    \ \ \mbox{ and }\ \      B_n (U_wQU_w) = B_n (Q) \circ \varphi_w
\eequ
for every $w \in \disc$.

The Berezin transform $B_0$ of operators was introduced by Berezin in \cite{ber} as tool to study spectral theory
and to construct approximations of the exponential of an operator. It has being used extensively to study properties
such as boundedness and compactness of Toeplitz, Hankel and other related operators.

The idea behind the transforms $B_n$ of functions in $L^1$ goes back to Berezin (see \cite{ber2}), and were explicitly used
in \cite{afr} to prove a deep result about the eigenfunctions of $B_0$ in the context of the ball in $\mathbb{C}^n$.
The extension of the definition of $B_n$ to operators is quite natural and appears in \cite{sua1}, where it is used to prove
approximation results in the same vein of Corollary \ref{aprochy} in the present paper.
\vspace{1mm}

The organization of the paper is as follows. In Section 2 we introduce the invariant Laplacian $\deltab$ and
prove some identities involving the interaction between $T_a^{(j)}$, $B_n$ and $\deltab$.
This will establish the technical foundations for the remaining sections.
In Section 3 we decompose $T_{B_n(S)}$ in terms of $T_{B_0(S)}^{(j)}$, and use it to give
a characterization of the $L^\infty$ closure of $B_0(\L(\berg))$, which turns out to be an algebra.
Section 4 contains the main results of the paper.
We prove that if $\mu \ge 0$ and $k\ge 1$, the operator $T_\mu^{(k)}$ is bounded (compact) if and only if $\mu$ is a
Carleson measure (resp.: a vanishing Carleson measure), and estimate the norms.
We also show that if $\mu$ is a complex measure whose variation $|\mu|$ is Carleson, then $T_\mu^{(k)}$ is the limit of
traditional Toeplitz operators. All these results generalize known facts for $k=0$.
In the last section we construct an example to show that for any $k\ge 0$,
$\|T_a^{(k+1)}\|$ is not majorized by $\sum_{j=0}^k \|T_a^{(j)}\|$  independently of $a\in L^\infty$.
In particular, the linear map $T_a \mapsto T_a^{(k+1)}$ is not bounded.
We will write indistinctly $T^{(0)}_a$ or $T_a$ for the traditional Toeplitz operator with symbol $a\in L^\infty$.

\section{The role of the invariant Laplacian}
\equnew

If $\Delta = \partial \ov{\partial}$ denotes a quarter of the usual Laplacian, where $\partial$ and $\ov{\partial}$
are the traditional Cauchy-Riemann operators, the invariant Laplacian is $\deltab := (1-|z|^2)^2 \Delta$.
It is easy to check that $(\deltab f)\circ \psi = \deltab (f\circ \psi)$ for every $f\in C^2(\disc)$ and
$\psi\in \mbox{Aut}(\disc)$. If $a\in L^\infty$ is such that $\deltab a\in L^1$, it is well known that
$\deltab B_0(a)= B_0(\deltab a)$. When also $\deltab a\in L^\infty$,  this equality rewrites as
$\deltab B_0(T_a)= B_0(T_{\deltab a})$. In accordance with this formula we give the following
\begin{fed}
Let
$$\D =\{ S\in \L(\berg):  \exists \, T\in \L(\berg) \,\mbox{ such that }\, \deltab B_0(S) = B_0(T) \} ,
$$
and define $\deltab: \D \, \rr \, \L(\berg)$ by\/ $\deltab S = T$.
\end{fed}
\noi This definition says that $\deltab B_0 (S) = B_0 (\deltab S)$ for all $S\in\D$. In \cite{sua1} it is showed that
if $S\in \L (\berg)$ and $n\geq 1$ then
\begin{equation}\laba{did1}
B_n (S) = \left( 1-\frac{ \tilde{\Delta} }{ n(n+1) } \right) B_{n-1} (S).
\end{equation}
Hence, a straightforward inductive argument shows that $\deltab B_n (S) = B_n (\deltab S)$ when $S\in\D$ for $n\ge 0$.
Also, the conformal invariance of $\deltab$ and \eqref{bobn} immediately prove that if $S\in \D$, then $U_wSU_w\in \D$ and
\bequ\laba{udeu}
\deltab (U_wSU_w) =  U_w (\deltab S) U_w.
\eequ
Observe also that \eqref{did1} implies that $\deltab B_n(S) \in L^\infty$ for every $S\in \L(\berg)$.

\begin{lemma}\laba{lew27}
Let $f,g,h,k$ be analytic on $\ov{\disc}$. Then
\begin{enumerate}
\item[{\em (i)}]  $\ \tilde{\Delta}(f\otimes g) = (f'\otimes g') + (z^2f)'\otimes (z^2g)' -2 \, (zf)'\otimes (zg)'$
\item[{\em (ii)}] $\ \la \tilde{\Delta} (f\otimes g) h,k \ra = \la \tilde{\Delta} (h\otimes k) f,g \ra .$
\end{enumerate}
\end{lemma}
\bdem [Proof of {\em (i)}]
\begin{align*}
\deltab B_0(f\otimes g)  &=  \deltab ( 1+|z|^4-2|z|^2) f \ov{g}
=  ( 1-|z|^2)^2 \big[  f'\ov{g'} + (z^2f)' \ov{(z^2g)'} -2 \, (zf)'\ov{(zg)'}  \big] \\*[1mm]
&=  B_0[ (f'\otimes g') + (z^2f)'\otimes (z^2g)' -2 \, (zf)'\otimes (zg)'  ] .
\end{align*}

\noi{\em Proof of {\em (ii)}.}
By (i),
\begin{equation}\label{dellato}
\tilde{\Delta} (z^n\otimes z^m) = nm(z^{n-1}\otimes z^{m-1})+ (n+2)(m+2) (z^{n+1}\otimes z^{m+1}) -2 (n+1)(m+1) (z^n\otimes z^m).
\end{equation}
Since $n  \|z^{n-1}\|^2 =1$ when $n>0$, for any $j, k \ge 0$ we have
\[  \la \tilde{\Delta} (z^n\otimes z^m)z^j, z^k\ra =
\left\{
\begin{array}{ll}
1     & \ \mbox{ if $(j,k) = (m-1,n-1)$} \\*[1mm]
\!-2  & \ \mbox{ if $(j,k) = (m,n)$} \\*[1mm]
1     & \ \mbox{ if $(j,k) = (m+1,n+1)$} \\*[1mm]
0     & \ \mbox{ otherwise}    \\*[1mm]
\end{array}
\right.
\]
This clearly shows that $\la \tilde{\Delta} (z^n\otimes z^m)z^j, z^k\ra =  \la \tilde{\Delta} (z^j\otimes z^k)z^n, z^m\ra$.
The lemma follows by sesqui-linearity.
\edem

\begin{lemma}\laba{cdover}
Let $\mu$ be a measure of finite variation such that $T_\mu^{(k)}$ is bounded for all $k\ge 0$.
Then $T_\mu^{(k)}\in\D$ for all $k\ge 0$, and
\bequ\laba{cfo}
\frac{\tilde{\Delta} T_\mu^{(k)}}{(k+1)} =  k \, T_\mu^{(k-1)}  + (k+2) \, T_\mu^{(k+1)} -2 (k+1) \, T_\mu^{(k)} ,
\eequ
or equivalently,  $(k+1)(k+2) \,  [T_\mu^{(k+1)} -  T_\mu^{(k)}] =
\deltab \Big[ T_{  \mu }^{(k)} + T_{  \mu }^{(k-1)}+ \cdots  + T_{  \mu }^{(0)} \Big]$.
Formally, we are taking $T_\mu^{(-1)}=0$ in \eqref{cfo} when $k=0$.
\end{lemma}
\bdem
By \eqref{dellato} with $k=n=m$,
\bequ\laba{deeka}
\frac{\tilde{\Delta} E_k}{k+1} =  k E_{k-1}+ (k+2) E_{k+1}-2 (k+1) E_k ,
\eequ
where $E_{-1}:=0$.
Since by \eqref{udeu}, $\tilde{\Delta} (U_w E_k U_w) = U_w (\tilde{\Delta}  E_k) U_w$, conjugating both members of
the above equality with respect to $U_w$ and integrating with respect to $(1-|w|^2)^{-2} d\mu(w)$, we obtain \eqref{cfo},
which is our claim.

The second formula follows from \eqref{cfo} by induction on $k$. It is immediate for $k=0$ and supposing
that it holds for an integer $k-1 \ge 0$, we get
\begin {align*}
\deltab T_{ \mu }^{(k)} + \deltab \Big[ T_{ \mu }^{(k-1)}+ \cdots  + T_{ \mu }^{(1)}+ T_{ \mu }^{(0)} \Big]
&= \deltab T_{  \mu }^{(k)} + k(k+1) [T_\mu^{(k)} -  T_\mu^{(k-1)}]\\
&= (k+1)(k+2) [T_\mu^{(k+1)} -  T_\mu^{(k)}] .
\end{align*}
Finally, if the last formula holds, substracting the equality for $k-1$ from the equality for $k$, we obtain
\eqref{cfo}.
\edem

\begin{lemma}
If $b_n, \, b\in L^\infty$ are such that $\|b_n\|_\infty \le C$, a constant independent of $n$,
and $b_n \to b$ pointwise, then $T_{b_n}^{(k)} \to T_{b}^{(k)}$ in the strong operator topology.
\end{lemma}
\bdem
We can assume that $b=0$. For $f, g\in \berg$,
\bequ\laba{mid}
|\la T_{b_n}^{(k)}f,g\ra|
\le \la T_{|b_n|}^{(k)} f, f \ra^{ \frac{1}{2} } \, \la T_{|b_n|}^{(k)} g, g \ra^{ \frac{1}{2} }
\le   \la T_{|b_n|}^{(k)} f, f \ra^{ \frac{1}{2} } \, C^{ \frac{1}{2} } \|g\|_2  ,
\eequ
where the first inequality follows from Cauchy-Schwarz's inequality and the second because
$\|T_{|b_n|}^{(k)}\| \le \|b_n\|_\infty \le C$. So, taking supremum in \eqref{mid} over $\|g\|_2=1$ for any fixed value of $n$,
we see that
$\| T_{b_n}^{(k)}f \|_2 \le C^{ \frac{1}{2} } \,  \la T_{|b_n|}^{(k)} f, f \ra^{ \frac{1}{2} } \to 0$ as $n\to\infty$
by the dominated convergence theorem.
\edem

\begin{propos}\laba{tk2}
Let $a\in L^\infty \cap C^2(\disc)$ such that $\deltab a\in L^\infty$. Then $\deltab T^{(k)}_a = T^{(k)}_{\deltab a}$.
\end{propos}
\bdem
For $0<r<1$ consider the functions $a_r(z)= a(rz)$. It follows from the previous lemma that
$T^{(j)}_{a_r} \to T^{(j)}_{a}$ in the strong operator topology when $r\to 1$ for all $j\ge 0$.
Then \eqref{cfo} implies that $\deltab T^{(k)}_{a_r} \tosot \deltab T^{(k)}_{a}$. Since
$(\deltab a_r) (z) = r^2  (\deltab a) (rz)$ is bounded by $\|\deltab a\|_\infty$, the previous lemma
says that $T^{(k)}_{\deltab a_r} \tosot T^{(k)}_{\deltab a}$.   Therefore it is enough to prove the lemma for $a_r$,
meaning that we can assume that $a\in C^2(\ov{\disc})$.  First observe that
$$
\deltab_z B_0(U_wE_kU_w)(z) =  \deltab B_0(E_k) (\varphi_w(z)) =  \deltab B_0(E_k) (\varphi_z(w))
= \deltab_w B_0(U_zE_kU_z)(w),
$$
where the equality in the middle holds because $\deltab B_0(E_k)$ is a radial function and
$|\varphi_w(z)| =  |\varphi_z(w)|$. Therefore
\begin{align*}
B_0(\deltab T^{(k)}_a )(w) &= \deltab B_0(T^{(k)}_a )(w)  = \int \deltab_w  B_0(U_zE_kU_z)(w) a(z) \, d\tilde{A}(z) \\
&=  \int \deltab_z  B_0(U_wE_kU_w)(z) a(z) \, d\tilde{A}(z) =  \int \Delta_z  B_0(U_wE_kU_w)(z) a(z) \, d{A}(z),
\end{align*}
and since $B_0(U_zE_kU_z)(w) =  B_0(U_wE_kU_w)(z)$ (because $B_0(E_k)$ is radial),
$$
B_0(T^{(k)}_{\deltab a} )(w)  =  \int   B_0(U_zE_kU_z)(w) \ (\deltab a)(z) \, d\tilde{A}(z)
=  \int   B_0(U_wE_kU_w)(z)\  (\Delta a)(z) \, d{A}(z) .
$$
Since for every fixed $w\in\disc$, the function
\bequ\laba{orp1}
B_0 (U_w E_k U_w) (z)= (1-|\varphi_w(z)|^2)^2 (k+1) (|\varphi_w(z)|^2)^k
\eequ
is defined for $z$ in some neighborhood of $\ov{\disc}$, the previous equalities and Green's theorem give
\begin{align*}
B_0\big( \deltab T^{(k)}_a - T^{(k)}_{\deltab a} \big) (w)
&= \int_\disc  \big[\Delta_z B_0 (U_w E_k U_w) (z) a(z) - B_0 (U_w E_k U_w) (z) (\Delta a)(z) \big] \, dA(z) \\
&= \int_{\partial \disc}
\left[  a(z)\, \frac{\partial}{\partial n}  B_0 (U_w E_k U_w) (z) - B_0 (U_w E_k U_w) (z) \, \frac{\partial a}{\partial n} (z)
\right] \frac{dm(z)}{\pi} ,
\end{align*}
where $\frac{\partial \ }{\partial n}$ is the derivative in the normal direction  and $dm(z)$ is the Lebesgue measure
on $\partial\disc$.  A straightforward calculation from \eqref{orp1} shows that both
$B_0 (U_w E_k U_w) (z)$ and $\frac{\partial}{\partial n}  B_0 (U_w E_k U_w) (z)$ vanish when $|z|=1$.
The Proposition follows from the fact that $B_0$ is one to one.
\edem

\begin{coro}\laba{harmostable}
If $a\in L^\infty$ is harmonic, $T_a^{(k)}= T_a$ for every integer $k\ge 0$.
\end{coro}
\bdem
By Proposition \ref{tk2}, $\deltab T_{ a }^{(k)}  =T_{ \deltab a }^{(k)}=0$ for all $k\ge 1$.
The corollary now follows from the second formula of Lemma \ref{cdover}.
\edem

\noi
Taking $a\equiv 1$ in the Corollary, we see that $T_1^{(k)}$ is the identity for all $k\ge 0$. This also follows from the
so called Schur orthogonality relations and it is the main ingredient in Englis's proof of the result
cited in the introduction. Indeed, the first inequality in \eqref{mid} implies that if $a\in L^\infty$, then
$\|T_a^{(k)}\| \le \|a\|_\infty \, \|T_1^{(k)}\|= \|a\|_\infty$.

\begin{propos}\laba{tk1}
Let $\mu$ be a finite measure such that $T^{(k)}_\mu$ is bounded for all $k\ge 0$.
Then $\ T_{ B_n(T_\mu^{(k)}) } = T^{(k)}_{ B_n(\mu) }$.
\end{propos}
\bdem
First we prove that $T_{ B_0(T_a^{(k)}) } = T^{(k)}_{ B_0(a) }$ by induction on $k$.
For $k=0$ there is nothing to prove. Suppose that the equality holds for $j=0, \ldots, k$.
By Proposition \ref{tk2}, the commutativity of $B_0$ and $\deltab$, and \eqref{cfo},
\begin{align*}
\deltab T_{ B_0(T_\mu^{(k)}) } &= T_{ \deltab B_0(T_\mu^{(k)}) }= T_{ B_0(\deltab T_\mu^{(k)}) }  \\
&=  (k+1)
\left[  k \, T_{  B_0 (T_\mu^{(k-1)}) } + (k+2) \, T_{  B_0 (T_\mu^{(k+1)})  } - 2 (k+1) \, T_{  B_0 (T_\mu^{(k)}) }  \right]
\end{align*}
and by \eqref{cfo},
$$
\tilde{\Delta} T_{B_0(\mu)}^{(k)}  =(k+1)
\left[ k \, T_{B_0(\mu)}^{(k-1)}  + (k+2) \, T^{(k+1)}_{  B_0 (\mu)  } -2 (k+1) \, T_{B_0(\mu)}^{(k)} \right] .
$$
By inductive hypothesis the left members of the above formulas are equal and we deduce that
$T_{  B_0 (T_\mu^{(k+1)})  }  = T^{(k+1)}_{  B_0 (\mu)  }$.

Now suppose that $k\ge 0$ is fixed and we prove the lemma by induction on $n$. So, suppose that the equality
holds for $n-1 \ge 0$. Then
\begin{align*}
n(n+1) [ T_{ B_{n-1} (T_\mu^{(k)}) }- T_{ B_{n}(T_\mu^{(k)}) }] &=
\deltab T_{ B_{n-1}( T_\mu^{(k)}) } =  \deltab T^{(k)}_{ B_{n-1}(\mu) }
=   n(n+1) [ T^{(k)}_{ B_{n-1}(\mu) }-T^{(k)}_{ B_{n} (\mu) }  ] .
\end{align*}
where the equality in middle holds by inductive hypothesis and the other two by Proposition \ref{tk2} and \eqref{did1}.
This proves our claim.
\edem

\section{$T_{B_n}$ in terms of $T^{(j)}_{B_0}$ and applications}

It is clear that $B_0: \L(\berg) \to L^\infty$ is not multiplicative but less clear that its image is not a multiplicative
set. We show this by constructing the following example.

Let $f,g \in\berg$ such that $T_f T_{\ov{g}}$ is bounded but $g\not \in \papa$. To see that such functions exist,
take for instance $f(z) = (1-z)^\alpha$ and $g(z) = (1-z)^{-\alpha}$, with $0<\alpha<1/2$. The elementary inequalities
$$
|1-z| \left(\frac{1-|w|}{1+|w|}\right)  \le  |1-\varphi_z(w)|   \le  |1-z| \left(\frac{1+|w|}{1-|w|}\right)
$$
yield
$$
B_0(|f|^p)B_0(|g|^p)(z) = \int |1-\varphi_z|^{p\alpha} dA \  \int \frac{dA}{ |1-\varphi_z|^{p\alpha} }
\le   \left[  \int \left(\frac{1+|w|}{1-|w|}\right)^{p\alpha} dA(w) \right]^2 <\infty
$$
if $0<p<\alpha^{-1}$. Hence, there is some $p>2$ such that $B_0(|f|^p)B_0(|g|^p)$ is bounded,
which by Theorem 5.2 of \cite{stzh} is a sufficient condition for the boundedness of $T_f T_{\ov{g}}$.

Since  $g\not \in \papa$, there is $h\in\berg$ such that $gh \not \in\berg$, implying that
the operator $(f\otimes gh)$ is not bounded. However, it is well defined on the reproducing kernels $K_z$  and satisfies
$(f\otimes gh) K_z =    \ov{g}(z) \ov{h}(z)  f \in \berg$ for all $z\in\disc$.
This holds because $K_z$ also reproduces functions in the Bergman space $A^1$.
In particular, its Berezin transform is defined,  and
$$
B_0(f\otimes gh)(z) = (1-|z|^2)^2 \ov{h}(z) \, f(z)\ov{g}(z)  = B_0(1\otimes h)(z) \, B_0(T_f T_{\ov{g}})(z).
$$
So, if $B_0(\L(\berg))$ is an algebra there must be $Q\in\L(\berg)$ such that $B_0(Q)=B_0(f\otimes gh)(z)$.
Consequently the function
$$F(z,w) := \la Q K_{\ov{z}} , K_{w} \ra - \la (f\otimes gh) K_{\ov{z}} , K_{w} \ra$$
is analytic on the bidisc $\disc^2$ and vanishes on the points $(z, \ov{z})$, implying that $F\equiv 0$.
Since the span of the reproducing kernels is dense in $\berg$, we conclude that $\|f\otimes gh\| = \|Q\| < \infty$,
a contradiction.

Despite the fact that $\L(\berg)$ is not an algebra, we will see that its closure is a uniform algebra, in fact,
the largest uniform algebra that previously known results allow. The key ingredient in the proof is the following
decomposition of $T_{B_n(S)}$, for $S\in\L(\berg)$.

\begin{lemma}\laba{decombn}
Let $S\in \L(\berg)$ and $n\ge 0$ integer. Then
\bequ\laba{bnbo1}
T_{B_n (S)} = (n+1) \sum_{j=0}^n {n \choose j}  \frac{(-1)^j}{j+1}  \    T^{(j)}_{ B_0(S) }.
\eequ
\end{lemma}
\bdem
\begin{align*}
B_0\left(  \sum_{j=0}^n {n \choose j}  \frac{(-1)^j}{j+1}  \,    T^{(j)}_{ B_0(S) } \right)\!(w)
&=      \!\int  \sum_{j=0}^n {n \choose j}  (-|\varphi_z(w)|^{2} )^j   (1-|\varphi_z(w)|^2)^2  B_0(S)(z) d\tilde{A}(z) \\
&=    \int \frac{ (1-|\varphi_z(w)|^2)^{n+2} }{ (1-|z|^2)^2 } \,  B_0(S)(z) dA(z) \\*[1mm]
&=  \frac{ B_n (B_0( S) )(w)}{(n+1) } \ = \ \frac{B_0( T_{B_n (S)} )(w)}{ (n+1) },
\end{align*}
where the last equality holds because $B_n$ and $B_0$ commute.
\edem

Consider the uniform algebra $\hip\subset L^\infty (\disc)$ of functions that are uniformly continuous from
the metric space $(\disc, \beta)$, where $\beta$ is the hyperbolic metric, into the complex plane
with the euclidean metric $(\mathbb{C}, |\,|)$.
In \cite{cob1} Coburn proved that $B_0(S)$ is a Lipschitz function between these metric spaces for every $S\in \L(\berg)$.
In particular, $B_0 (\L(\berg)) \subset \hip$, a fact used in \cite{sua1} to study some subalgebras of
$\L(\berg)$ in terms of their Berezin transforms. We see next that the inclusion is dense.

\begin{theo}
The $L^\infty$-closure of $B_0(\L(\berg))$ is $\hip$.
\end{theo}
\bdem
Let  $a\in \hip$. Replacing $B_0(S)$ by $a$ in the chain of equalities of the previous proof (except for the last one), gives
$$
B_0\left( (n+1) \sum_{j=0}^n {n \choose j}  \frac{(-1)^j}{j+1}  \    T^{(j)}_{ a } \right)  = B_n(a) .
$$
Taking $d\mu = adA$ in  \eqref{dmuda}, a change of variables shows that

$$
B_n(a)(z)  =  \int_D a( \varphi_z (\zeta ) )   (n+1) (1-|\zeta |^2 )^n \, dA(\zeta ) \to a(\varphi_z(0))=a(z).
$$
uniformly on $z$ when $n\to \infty$, because since $a\in\hip$, the functions $a\circ\varphi_z$ are equicontinuous at $0$,
and the probability measures $(n+1)(1-| . |^2)^n dA$ tend to accumulate all the mass at $0$ when $n\to\infty$.
Thus, $\hip \subset \ov{B_0(\L(\berg))}$.
\edem

\begin{coro} The set\/ $\{ T_{B_0(S)} : \  S\in  \L(\berg)  \}$    is norm dense in $\{ T_a : \  a\in L^\infty \}$.
\end{coro}
\bdem
The last corollary implies that the first set is norm dense in $\{ T_a : \  a\in \hip \}$,
which by \cite[Thm.$\,$5.7]{sua1} is norm dense in $\{ T_a : \  a\in L^\infty \}$.
\edem

The next result is an easy consequence of the identities in the previous section and Lemma \ref{decombn}.
We need some notation first. Let $m\geq 0$ be an integer and\/ $x=\{ x_n \}_{n\geq 0}$ be a sequence of complex numbers.
The \mbox{$m$-difference} of $x$, denoted $\Delta^m x$, is the sequence whose $n$-th term is
$$
\Delta^m_n x:= (-1)^m \sum_{j=0}^m   {m \choose j}  (-1)^j \,x_{n+j} , \ \mbox{ for }\   n\geq 0 .
$$
That is, $\Delta^m$ is the $m$-iteration of the difference operator $\Delta \{ x_n \}_{n\geq 0} :=\{ x_{n+1}- x_n \}_{n\geq 0}$.

\begin{propos} Let $f, g,  h,  k\in\berg$ and integers $n, j\ge 0$. Then
$$\la T_{B_n (f\otimes g)}h,k\ra  =  \la   T_{ B_n (h\otimes k)} f , g\ra$$
$$
\hspace{-10mm}\mbox{and} \hspace{6mm} \int    \la U_w e_j , h \ra   \ov{\la  U_w e_j , k\ra}    \   f(w)\ov{g(w)} dA(w)
=  \int   \la U_w e_j , f \ra   \ov{\la  U_w e_j , g\ra}    \   h(w)\ov{k(w)} dA(w).
$$
\bequ\label{miacle}
\hspace{0mm}\text{In particular,}\hspace{6mm}
\int |\la U_w e_j , h \ra|^2   \  |f(w)|^2 dA(w) = \int |\la U_w e_j , f \ra|^2   \  |h(w)|^2 dA(w) .
\eequ
\end{propos}
\bdem
Since $\|T_{B_n (f\otimes g)}\| \le  C_n \|f\|_2 \, \|g\|_2$, it is enough to assume that all the functions are polynomials.
Since $B_0 (f\otimes g)=(1-|z|^2)^2 f\ov{g}$, the first assertion is clear for $n=0$.
So, assuming that the result holds up to $n$, by \eqref{did1} we need to prove the equality for $\deltab B_n$ instead of $B_n$.
\begin{align*}
\la \deltab T_{B_n (f\otimes g)}h,k\ra
&=  \la  B_n (f\otimes g) h,k\ra  +  \la  B_n ((z^2f)'\otimes (z^2g)') h,k\ra   -2 \la  B_n ((zf)'\otimes (zg)') h,k\ra\\
&=  \la  B_n (h\otimes k) f,g\ra  +  \la  B_n (h\otimes k)(z^2f)' ,(z^2g)'\ra   -2 \la  B_n (h\otimes k) (zf)',(zg)'\ra\\
&=  \la  B_n (h\otimes k) , \Delta (1-|z|^2)^2 \ov{f}g\ra \\
&=  \la  \deltab B_n (h\otimes k) ,  \ov{f}g\ra  \\
&=  \la  \deltab T_{ B_n (h\otimes k)} f , g\ra ,
\end{align*}
where the first equality follows from Proposition \ref{tk2}, the commutativity of $B_n$ and $\deltab$,
and
Lemma \ref{lew27}, the second equality holds by inductive hypothesis, the fourth one by Green's theorem,
and the last one by Proposition \ref{tk2} again. Writing $\sigma_j(z) = z^j$, \eqref{bnbo1} says that
\begin{align*}
\frac{ T_{B_n (f\otimes g)} }{n+1}&=  \sum_{j=0}^n {n \choose j}  \frac{(-1)^j}{j+1}  T^{(j)}_{ B_0(f\otimes g) }
= (-1)^n  \int      \Delta_0^n  (U_w\sigma_j\otimes U_w\sigma_j)    \  f(w)\ov{g(w)} dA(w).
\end{align*}
Therefore the equality $\la T_{B_n (f\otimes g)}h,k\ra =  \la T_{B_n (h\otimes k)}f,g\ra$ rewrites as
$$
\Delta^n_0 \int  \la U_w \sigma_j , h \ra   \ov{\la  U_w \sigma_j , k\ra}   \,   f(w)\ov{g}(w) dA(w) =
\Delta^n_0\int   \la U_w \sigma_j , f \ra   \ov{\la  U_w \sigma_j , g\ra}    \,   h(w)\ov{k}(w) dA(w),
$$
and the second claim follows by induction on $n$.
\edem

\section{Carleson measures as symbols}
\equnew

A positive measure $\mu$ on $\disc$ is called a Carleson measure if $\berg \subset L^2(d\mu)$.
If in addition the inclusion is compact, $\mu$ is called a vanishing Carleson measure.
Among the many known characterizations of Carleson measures (see \cite[p.$\,$123]{zhu} for comments and references),
a positive measure $\mu$ is Carleson if and only if $\| B_0(\mu) \|_\infty < \infty$, a quantity that is equivalent to
the operator norm of the inclusion of $\berg$ in $L^2(\mu)$.
Another characterization comes from replacing the kernel of the Berezin integral by a box kernel.
Indeed, if $0\le r<1$ and $v\in\disc$, consider the pseudo-hyperbolic disk
$$
D(v,r):= \{ z\in\disc: \, |\varphi_v(z)| \le r  \}\ \ \mbox{  and its area $\ \ |D(v,r)|:= \int_{D(v,r)} dA.$}$$
If $\mu$ be a positive measure on $\disc$ and $0<r<1$, there is a constant $C(r)>0$ depending only on $r$ such that
\bequ\laba{triui}
\frac{1}{C(r)} \sup_{v\in\disc}\frac{\mu(D(v,r)) }{|D(v,r)|}\leq    \| B_0(\mu)\|_\infty
\leq  C(r)\sup_{v\in\disc}\frac{\mu(D(v,r)) }{|D(v,r)|}.
\eequ
Clearly, if the above  supremum is finite for some $r$ then it is finite for all $1<r<1$.
Finally, a positive measure $\mu$ is Carleson if and only if $T_\mu$ is bounded
(see \cite[pp.$\,$111-112]{zhu}). We shall see that the same holds for  $T_\mu^{(k)}$ when $k\ge 1$.
For a positive measure $\mu$  write $d\tilde{\mu} := (1-|z|^2)^{-2} d\mu$.
\vspace{2mm}

\begin{lemma}\laba{potriui}
Let $\mu$ be a positive finite measure on $\disc$. Then
$$
\frac{\mu(D(v,r)) }{|D(v,r)|}    \left[ \frac{r(1-r^2) }{ 4 } \right]^2
\le \tilde{\mu}(D(v,r)) \le   \frac{\mu(D(v,r)) }{|D(v,r)|}    \left[ \frac{4r }{(1-r^2)^2}  \right]^2
$$
for every $v\in\disc$ and $0<r<1$.
\end{lemma}
\bdem
Since by \cite[p.$\,$60]{zhu}, $|D(v,r)|  = \left[\frac{ r(1-|v|^2)  }{ 1-|v|^2 r^2 }\right]^2$,
\begin{align*}
\tilde{\mu}(D(v,r)) &= \int_{D(v, r)} \frac{d\mu(\xi)}{(1-|\xi|^2)^2}
= \frac{1}{|D(v,r)|}  \int_{D(v, r)}  \left[ \frac{r(1-|v|^2) }{(1-|\xi|^2) (1-|v|^2 r^2)} \right]^2  d\mu(\xi).
\end{align*}
The lemma follows immediately from the easy inequalities, valid for $\xi\in D(v,r)$:
$$
\frac{(1-r^2) }{ 4 } \le  \frac{(1-|v|^2) }{(1-|\xi|^2) }  \le \frac{4 }{(1-r^2)}  . \vspace{-5mm}
$$
\edem

\begin{theo}\laba{Carmer}
Let $\mu$ be a positive finite measure on $\disc$. Then $T_\mu^{(k)}$ is bounded if and only if
$\mu$ is a Carleson measure, in which case,
\bequ\laba{atknor}
\frac{C}{(k+2)} \ \| B_0(\mu)\|_\infty  \leq \| T^{(k)}_\mu \| \leq  4(k+2)\| B_0(\mu) \|_\infty ,
\eequ
where $C>0$ is an absolute constant.
\end{theo}
\bdem
First let us assume that $T^{(k)}_\mu$ is a bounded operator. For $k\geq 1$ consider the function $f(x)= (k+1) x^k (1-x)^2$
defined in $[0,1]$. This function reaches its maximum at $x= k/(k+2)$.
If $\,\frac{k-1/2}{k+2}\leq x\leq \frac{k+1}{k+2}$ (that is, $x= \frac{k+y}{k+2}\,$ with $-1/2 \leq y\leq 1$), then
$$
f(x)= f(\frac{k+y}{k+2})= (k+1) \left[\frac{k+y}{k+2}\right]^k \left[\frac{2-y}{k+2}\right]^2
\ge  \frac{(k+1)}{(k+2)^2}  \left[1-\frac{5/2}{k+2}\right]^k \ge  \frac{ c_1 }{(k+2)},
$$
where $c_1 >0$ is a constant independent of $k$.
This means that there is an absolute constant $c_1 >0$ such that for all $k\ge 1$,
\bequ\laba{qp1}
  (k+1) |z|^{2k} (1-|z|^2)^2 \ge  \frac{c_1}{(k+2)}
 \ \ \ \mbox{ if }\ \ \  \frac{k-1/2}{k+2}\leq |z|^2 \leq \frac{k+1}{k+2}.
\eequ
Now, let  $0< r\leq z_k := \sqrt{\frac{k}{k+2}}$. By the geometrical arguments in \cite[p.$\,$3]{gar},
$D(z_k, r)$ is contained in the annulus
$$
\frac{z_k -r}{1-rz_k} \le |w| \le \frac{z_k +r}{1+rz_k}.$$
Thus, if we choose $r \leq \sqrt{\frac{k}{k+2}}$ small enough so that
\bequ\laba{srtq}
\sqrt{\frac{k-1/2}{k+2}} \le  \frac{  \sqrt{\frac{k}{k+2}} -r }{1-r \sqrt{\frac{k}{k+2}} }
\ \ \ \mbox{ and }\ \ \
\frac{ \sqrt{\frac{k}{k+2}} +r}{1+r  \sqrt{\frac{k}{k+2}} } \le  \sqrt{\frac{k+1}{k+2}}
\eequ
for all $k\ge 1$, then  $D(z_k, r)$ is contained in the annulus $\frac{k-1/2}{k+2}\leq |z|^2 \leq \frac{k+1}{k+2}$,
implying that the inequalities in \eqref{qp1} hold for $z\in D(z_k, r)$.
We see next that $0<r\leq 1/10$ does the trick. Clearing $r$ from \eqref{srtq} we get the equivalent inequalities
$$
r   \le  \frac{ \sqrt{\frac{k}{k+2}}-\sqrt{ \frac{k-1/2}{k+2} } }{[1 -  \sqrt{\frac{k}{k+2}} \sqrt{ \frac{k-1/2}{k+2} }]}
\ \ \mbox{ and }\ \
r   \le  \frac{ \sqrt{ \frac{k+1}{k+2} }- \sqrt{\frac{k}{k+2}} }{[1-  \sqrt{\frac{k}{k+2}}\sqrt{ \frac{k+1}{k+2} }]},
$$
or equivalently,
$$
r \le \min \left\{  \frac{ \sqrt{k+2} }{[\sqrt{k}+\sqrt{ k-1/2 }]} \, \frac{1/2 }{[k+2 -  \sqrt{ k^2-k/2 }]}, \,
\frac{ \sqrt{k+2} }{[\sqrt{k}+\sqrt{ k+1 }]} \, \frac{ 1 }{[k+2- \sqrt{ k^2+k }]}  \right\} .
$$
The claim follows because this minimum is bounded below by
\begin{align*}
\frac{ \sqrt{k+2} }{[\sqrt{k}+\sqrt{ k+1 }]} \,  \frac{1/2 }{[k+2 -  \sqrt{ k^2-k/2 }]}
&\ge  \frac{1/4 }{[k+2 -  \sqrt{ k^2-k/2 }]}   \ = \ \frac{k+2 +  \sqrt{ k^2-k/2 } }{ [18k+16]}\\*[1mm]
&\ge  \frac{2k+3/2   }{ [18k+16]} \ \ge\  \frac{2+3/2   }{ [18+16]} \ > \ \frac{1}{10}.
\end{align*}
Therefore
\begin{align}\laba{jjt}
B_0(T^{(k)}_\mu) (w)
&\stackrel{\mbox{\scriptsize{\hspace{12mm}}}}{=}
\int (k+1) |\varphi_w(z)|^{2k} (1-|\varphi_w(z)|^2)^2 \frac{d\mu(z)}{(1-|z|^2)^2} \nonumber \\
&\stackrel{\mbox{\scriptsize{\hspace{12mm}}}}{\ge}  \int_{ D (\varphi_w(z_k) , r)  } (k+1)
|\varphi_w(z)|^{2k} (1-|\varphi_w(z)|^2)^2 \frac{d\mu(z)}{(1-|z|^2)^2} \nonumber \\
&\stackrel{\mbox{\scriptsize{by \eqref{qp1}}}}{\ge}
\frac{c_1}{(k+2)}  \int_{ D (\varphi_w(z_k) , r)  }  \frac{d\mu(z)}{(1-|z|^2)^2} \nonumber \\
&\stackrel{\mbox{\scriptsize{\hspace{12mm}}}}{=}   \frac{c_1}{(k+2)}  \tilde{\mu}( D(  \varphi_w( z_k  ) , r)  ).
\end{align}
Taking  the supremum for $w\in \disc$ and using that  $\{\varphi_w( z_k ) : \ w\in\disc\} =\disc$
for any fixed $z_k\in\disc$, we get
\bequ\laba{caqui}
\| T^{(k)}_\mu \| \ge \| B_0(T^{(k)}_\mu) \|_\infty \ge
\frac{c_1}{(k+2)} \sup_v \tilde{\mu}\left( D \left(  v , r\right)  \right)
\eequ
for any $r\le 1/10$.

By \eqref{triui}, Lemma \ref{potriui} and \eqref{caqui}, there are absolute constants $C_0, \, C_1$ and $C_2$, such that
$$
\| B_0(\mu)\|_\infty \le C_0  \sup_v \frac{\mu(D(v,\frac{1}{10} )) }{|D(v,\frac{1}{10} )|}
\le   C_1 \sup_v  \tilde{\mu}(D(v,\frac{1}{10} ))  \le C_2 (k+2) \| T^{(k)}_\mu \| .
$$
This proves the first inequality in \eqref{atknor}. \\

\noi
Now suppose that $\mu$ is a Carleson measure, and let $F(z) = \sum a_j e_j(z) \in \berg$.
For $0\le t < 2\pi$ and $0\le r<1$ we have
\begin{eqnarray*}
|\la F(e^{it}z) , (U_re_k)(z)\ra|^2    &=&
\sum_{j,l}   a_j \ov{a}_l  \,  \la e_j(e^{it}z) , (U_re_k)(z) \ra  \,  \ov{\la e_l(e^{it}z) , (U_re_k)(z) \ra} \\
&=&  \sum_{j,l}   a_j \ov{a}_l  \, e^{i(j-l)t}\,  \la e_j, U_re_k \ra   \,  \ov{\la e_l, U_re_k\ra} ,
\end{eqnarray*}
and since
$\,  |\la F, U_{re^{it}}e_k\ra|  =  |\la F(z), (U_re_k)(e^{-it}z)\ra|  =  |\la F(e^{it}z), (U_re_k)(z)\ra|$,
then
\begin{eqnarray*}
\int_0^{2\pi}  |\la F, U_{re^{it}}  e_k\ra|^2 \,\frac{dt}{2\pi}
&=&   \sum_{j}   |a_j|^2  |\la e_j, U_re_k \ra|^2   \ \geq\        |a_k|^2  |\la e_k, U_re_k \ra|^2 \\
&=&  |\la F, e_k\ra|^2  |\la e_k, U_re_k \ra|^2    \\*[1mm]
&=&  |\la F, e_k\ra|^2    \int_0^{2\pi}  |\la e_k, U_{re^{it}}e_k\ra|^2 \,\frac{dt}{2\pi} .
\end{eqnarray*}
Multiplying by $2r dr$ and integrating yields
$$
\int  |\la F, U_ze_k\ra|^2 \,dA(z)  \geq  |\la F, e_k\ra|^2    \int   |\la e_k, U_ze_k\ra|^2\, dA(z) .
$$
So, taking $F= U_wf$ we get
$$
\int  |\la U_wf, U_ze_k\ra|^2 \, dA(z) \geq  |\la U_wf, e_k\ra|^2    \int |\la e_k, U_ze_k\ra|^2\,dA(z) .
$$
Writing $\lambda =  (z\ov{w} -1) / (1-w\ov{z} )$, we have $U_w U_z= U_{\varphi_{w}(z)}  V_\lambda $,
where $(V_\lambda h) (\om) = \lambda h(\lambda \om)$ for $h\in \berg$.  Consequently,
$|\la U_wf, U_ze_k\ra| = |\la f, U_wU_ze_k\ra| = |\la f, U_{\varphi_w(z)}  e_k\ra|$,
and the change of variables $v= \varphi_w(z)$ in the first integral above yields
$$
\int  |\la f, U_ve_k\ra|^2 \, |\varphi'_w(v)|^2   dA(v)  \geq
|\la U_wf, e_k\ra|^2  \,  \int   |\la e_k, U_ze_k\ra|^2\, dA(z) .
$$
Integrating with respect to $d\tilde{\mu}(w)$,
\bequ\laba{iwrt}
 \int_{\disc} \left[\int \frac{(1-|v|^2)^2}{|1-\ov{w}v|^4} d\mu(w)\right] \ |\la f, U_ve_k\ra|^2 \, d\tilde{A}(v)
 \geq  c_k\int_\disc |\la U_wf, e_k\ra|^2 d\tilde{\mu}(w),
\eequ
where
\begin{align*}
c_k &= \int   |\la e_k, U_ze_k\ra|^2\, dA(z) \!\stackrel{\mbox{\scriptsize{by \eqref{miacle}}}}{=}\!
\int |\la U_z e_k , 1 \ra|^2   \,  |e_k(z)|^2 dA(z) =  \int (1-|z|^2)^2   \,  |e_k(z)|^4 dA(z)\\
&=  (k+1)^2 \int_0^1  (1-x)^2   \,  x^{2k} dx  = \frac{ (k+1)^2 \, 2! \, (2k)!}{(2k+3)!}
= \frac{ (k+1)   }{(2k+3)(2k+1)} \ge \frac{ 1  }{4(k+2)},
\end{align*}
where the solution to the integral comes from $\int_0^1 (1-t)^p \, t^q \,dt =  p! \, q!/(p+q+1)!$ for integers $p, q\ge 0$.
Thus, going back to \eqref{iwrt},
\begin{align*}
\| B_0(\mu) \|_\infty\   \|f\|^2 &\ge \la T_{B_0(\mu)}^{(k)}f, f\ra
\geq   c_k \la T_\mu^{(k)} f,f \ra
\geq      \frac{ 1  }{4(k+2)} \ \la T_\mu^{(k)} f,f \ra     .
\end{align*}
This proves the second inequality in \eqref{atknor}.
\edem

\noi
It would be interesting to know how sharp are the bounds in \eqref{atknor} except for absolute multiplicative constants when $k$
tends to infinity, especially the upper bound.

\begin{rema}\laba{remate}
Observe that by \eqref{caqui} and the subsequent inequality, we also showed that
$$
\frac{C}{(k+2)} \| B_0(\mu)\|_\infty \le  \| B_0(T^{(k)}_\mu) \|_\infty \le  \| T^{(k)}_\mu \|  ,$$
and that the last formula of the proof says that\/ $4(k+2)T_{B_0(\mu)}^{(k)} \ge  T_\mu^{(k)}$ as positive operators.
\end{rema}

\vspace{1mm}

\noi
Suppose that $\mu$ is a complex measure on $\disc$ such that its variation $|\mu|$ is Carleson.
By \eqref{mid} with measures instead of functions, we see that $\| T^{(k)}_\mu \| \leq \| T^{(k)}_{|\mu|} \|$ for all $k\ge 0$,
so  $T^{(k)}_\mu\in \L(\berg)$.
It is worth noticing that the converse does not hold, since there are finite measures $\mu$
such that  $T_\mu$ is bounded but  $|\mu|$ is not Carleson.
The next result was proved in \cite[Cor.$\,$2.5]{sua3} for $k=0$.   In particular, it shows that
when $a\in L^\infty$, $T_a^{(k)}$ is a limit of classical Toeplitz operators.

\begin{coro}\laba{aprochy}
Let $\mu$ be a finite measure on $\disc$ such that $|\mu|$ is a Carleson measure and $k\ge 0$ be an integer. Then
$$
T_{ B_n( T_\mu^{(k)} ) } \rr T_\mu^{(k)} \ \ \mbox{ when $\, n\rr \infty$.}
$$
\end{coro}
\bdem
Decomposing $\mu= \mu_1 + i \mu_2$, where each $\mu_j$ is a real measure, and using Jordan decomposition with both
$\mu_1$ and $\mu_2$,
we can assume without loss of generality that $\mu \ge 0$.
By Lemma 4.1 of \cite{sua2}, if  $Q\in \L(\berg)$ satisfies $\|T_{\deltab B_n(Q)}\| \leq\/ C$, where $C$ is
independent of\/ $n$, then $T_{B_n(Q)} \rr Q$. So, we need to prove the above inequality for $Q=T_\mu^{(k)}$.

By Propositions \ref{tk2} and \ref{tk1},  and \eqref{cfo},
$$
T_{ \deltab B_n( T_\mu^{(k)}) } = \deltab  T_{ B_n( T_\mu^{(k)}) } = \deltab  T^{(k)}_{ B_n(\mu) }
=(k+1)  [ k \, T_{ B_n(\mu) }^{(k-1)}  + (k+2) \, T_{ B_n(\mu) }^{(k+1)} -2 (k+1) \, T_{ B_n(\mu) }^{(k)} ].
$$
Since $B_n(\mu) dA$ is a Carleson measure satisfying
$\| B_0 B_n (\mu) \|_\infty = \| B_n B_0 (\mu) \|_\infty \le \| B_0(\mu) \|_\infty$,
using \eqref{atknor} in the above equality gives
$\| T_{ \deltab B_n( T_\mu^{(k)}) } \|  \le  4^2  (k+3)^3 \| B_0(\mu) \|_\infty$, which does depend on $n$.
\edem

It is well known that a positive measure $\mu$ on $\disc$, the condition of being a vanishing Carleson measure
is equivalent to $B_0(\mu)(z)\to 0$ when $|z|\to 1$, and also to the compactness of $T_\mu$
(see \cite[pp.$\,$112-115]{zhu},  also \cite[Propo.$\,$3]{McS}).
We aim to prove the same result for $T^{(k)}_\mu$ when $k$ is any nonnegative number.

\begin{lemma}
If $f_n\in \berg$ is a sequence that tends weakly to $0$ then $\la f_n , U_w e_k \ra \rr 0$
uniformly for $w$ in compact sets of\/ $\disc$.
\end{lemma}
\bdem
By the Banach-Steinhaus Theorem (see \cite[p.$\,$44]{rud}) the norms $\|f_n\|$ are uniformly bounded
and by Lemma 4.3 of  \cite{sua1} the function $w\mapsto U_w e_k$ is uniformly continuous on compact sets.
Thus, the Cauchy-Schwarz inequality shows that the scalar functions $F_n(w) = \la f_n , U_w e_k \ra$
are equicontinuous on compact sets. Since by hypothesis $F_n \rr 0$ pointwise,
Ascoli's theorem (see \cite[p.$\,$394]{rud}) implies that $F_n \rr 0$ uniformly on
compact sets.
\edem

\begin{lemma}\laba{preh}
If $a\in L^\infty$ has compact support then $T_a^{(k)}$ is compact.
\end{lemma}
\bdem
Let $f_n\in\berg$ be a sequence that tends weakly to $0$. Then
$$
|\la T_a^{(k)}f_n, f_n\ra| \le
\|a\|_\infty \, \tilde{A}(\supp a)  \sup_{w\in {\scriptsize \supp} a} |\la f_n , U_w e_k \ra|^2,
$$
whose last factor tends to $0$ by the previous lemma.
\edem

\begin{theo}
Let $\mu$ be a positive finite measure on $\disc$. Then $T_\mu^{(k)}$ is compact if and only if
$\mu$ is a vanishing Carleson measure.
\end{theo}
\bdem
Suppose that $\mu$ is a vanishing Carleson measure and let $0<r<1$. By Remark \ref{remate},
$$
0\le T_\mu^{(k)}  \le 4(k+2) T^{(k)}_{ B_0(\mu) }  = 4(k+2) \left[ T^{(k)}_{ \chi_{rD}B_0(\mu) } +
T^{(k)}_{ \chi_{D\setminus rD}B_0(\mu) }\right].
$$
By Lemma \ref{preh} the first operator in the sum is compact and by Englis's theorem,
$$
\| T^{(k)}_{ \chi_{D\setminus rD}B_0(\mu) } \| \le \| \chi_{D\setminus rD}B_0(\mu)\|_\infty \rr 0
\ \mbox{ when $r\rr 1$}.
$$
Thus, $T_\mu^{(k)}$ is compact. Conversely, suppose now that $T_\mu^{(k)}$ is compact.
Then $B_0(T^{(k)}_\mu) (w)\rr 0$ when $|w|\rr 1$, which together with \eqref{jjt} says that
there are $z_k \in \disc$ and $0<r<1$ such that
$$\tilde{\mu}( D(  \varphi_w( z_k  ) , r)  ) \rr 1 \mbox{ when $|w|\rr 1$}.
$$
If $V\subset \disc$ is such that $\disc\setminus V$ is compact, the same holds for the set
$\{\varphi_w( z_k ) : \ w\in V \}$, for any fixed $z_k\in\disc$. Therefore
$\tilde{\mu}( D(  v , r)  ) \rr 1$  when $|v|\rr 1$, which together with Lemma \ref{potriui}
gives
$$\frac{\mu(D(v,r)) }{|D(v,r)|}    \to 0 \ \mbox{ as }\   |v|\to 1 .$$
Then $\mu$ is a vanishing Carleson measure by \cite[pp.$\,$111-114]{zhu} .
\edem

\section{Example of bad behaviour}
\equnew

As far as I know there is no accurate estimate for $\|T_a\|$ when $a\in L^\infty$ is arbitrary,
which obviously remains true for $\|T^{(k)}_a\|$ when $k\ge 1$.
It would be interesting to know if at least $\|T^{(k)}_a\|$ is majorized by $\|T_a\|$, or more generally,
if for some given $k\ge 1$, there exists a positive constant $C_k$ depending
only on $k$ such that
\bequ\laba{falsy}
\| T_a^{(k)} \| \le C_k ( \| T_a^{(0)} \| + \cdots + \| T_a^{(k-1)} \|  )
\ \mbox{ for all $a\in L^\infty$}.
\eequ
By Theorem \ref{Carmer} this is certainly the case when $a\ge 0$ or when $adA$ is replaced by any Carleson measure.
Unfortunately \eqref{falsy} does not hold for any $k\ge 1$, as the example that we construct next will show.

\begin{lemma}\laba{lincom}
For $a\in L^\infty$ and $\ell \ge 0$ there are constants $c_0, \ldots , c_{\ell}$ depending only on $\ell$ such that
$$T_a^{(\ell)}= c_0 \deltab^0 T_a + \cdots + c_{\ell} \deltab^{\ell} T_a .
$$
\end{lemma}
\bdem
By the second formula of Lemma \ref{cdover},
\begin{align*}
T_a^{(\ell)}  &= T_a^{(0)}  +
\deltab \, \sum_{m=0}^{\ell-1} \frac{1}{(m+1)(m+2)} \Big[ T_{  a }^{(m)} + T_{  a }^{(m-1)}+ \cdots  + T_{  a }^{(0)} \Big].
\end{align*}
This proves the lemma for $\ell=1$  and assuming inductively that it holds for $T_a^{(m)}$ with $m=1, \ldots, \ell-1$,
it also shows that it holds for $T_a^{(\ell)}$.
\edem
\begin{coro}
For all $k\ge 0$ and $a\in L^\infty$ there is $C_k>0$ such that
$$\sum_{\ell=0}^{k} \| T^{(\ell)}_a\| \le C_k \sum_{\ell=0}^{k} \| \deltab^{\ell} T_a\| .  $$
\end{coro}

\noi
The proof of Lemma \ref{lincom} clearly shows that both the lemma and its corollary hold if $adA$ is replaced by any
finite measure $\mu$ such that $T^{(k)}_\mu$ is bounded for every $k\ge 0$. In particular, they hold when $|\mu|$
is a Carleson measure.

Let $k\ge 1$ and suppose that \eqref{falsy} holds. This, together with \eqref{cfo} imply the first of the following
inequalities
$$
\| \deltab^k T_a \| \le  C_1(k) \sum_{\ell=0}^{k-1} \| T^{(\ell)}_a\| \le C_2(k) \sum_{\ell=0}^{k-1} \| \deltab^{\ell} T_a\|
\ \mbox{ for all $a\in L^\infty$,}$$
for some $C_1(k) >0$, where the second inequality comes from the corollary. Thus, the next example disproves \eqref{falsy}.

\vspace{1mm}
\noi {\bf Example.} We claim that if $k\ge 1$ there is no positive constant $C_k$ such that
 $$\| \deltab^k T_a \| \le C_k \sum_{\ell=0}^{k-1} \| \deltab^{\ell} T_a\|
 \hspace{5mm}   \mbox{ for all } a\in L^\infty .
 $$
For $j\ge 0$ recall that $E_j = e_j\otimes e_j$, and we write $E_j=0$ if $j<0$. An iteration of \eqref{deeka} shows that
$\tilde{\Delta}^\ell E_j$ is a linear combination of $E_{j-\ell}, \ldots , E_{j+\ell}$ in such a way that
there are positive constants $c_\ell$ and $C_\ell$ independent of $j$ with
$c_\ell (j+1)^{2\ell} \le \| \tilde{\Delta}^\ell E_j \| \le C_\ell (j+1)^{2\ell}$ for all $\ell \ge 0$. In particular,
if $0\le \ell \le k$, there are constants $c$ and $C$ depending only on $k$ such that
$$c (j+1)^{2\ell} \le \| \tilde{\Delta}^\ell E_j \| \le C (j+1)^{2\ell}    \hspace{1cm} \forall \ell=0,\ldots, k
\,\mbox{ and }\, j\ge 0.
$$
By \cite[Thm.$\,$4.3]{sua2}, $T_{B_n(E_j) } \to E_j$ when $n\to\infty$.
Hence, Proposition \ref{tk2}, the commutativity of $B_n$ and $\deltab$, and the previous comments  yield
$$
\tilde{\Delta}^\ell T_{B_n(E_j) }  =  T_{\tilde{\Delta}^\ell  B_n(E_j) } =  T_{B_n( \tilde{\Delta}^\ell  E_j) }
\to \tilde{\Delta}^\ell  E_j, \ \ \mbox{ as } n\to \infty .
$$
Therefore for each pair of integers $k, j \ge 0$ we can choose $n=n(k, j)$ large enough so that
$$
\frac{c}{2} (j+1)^{2\ell} \le \| \tilde{\Delta}^\ell T_{B_n(E_j) } \| \le 2C (j+1)^{2\ell}  \hspace{1cm}
\forall \ell=0,\ldots, k.
$$
Taking $a_j:= (j+1)^{-2k} B_n(E_j)\in L^\infty$, the above inequalities show that,
$$
\sum_{\ell=0}^{k-1} \| \tilde{\Delta}^\ell T_{a_j} \| \le 2C \sum_{p=1}^{k} \frac{1}{(j+1)^{2p}}
\le    \frac{2C}{(j+1)^2-1},
\ \mbox{ while }\
\frac{c}{2} \le \| \tilde{\Delta}^k T_{a_j} \|
$$
for all $j\ge 1$. Taking $j\to \infty$ shows our claim.\\

\noindent {\bf Acknowledgement:}
Research supported in part by the ANPCyT grant PICT2009-0082 and
UBA grant UBACyT 20020100100502, Argentina.

\newcommand{\foo}{\footnotesize}

 \noindent Daniel Su\'{a}rez\\
 Departamento de Matem\'{a}tica \\
 Facultad de Cs.\@ Exactas y Naturales \\
 UBA, Pab. I, Ciudad Universitaria \\
 (1428) N\'{u}\~{n}ez, Capital Federal \\
 Argentina\\
\vspace{0.5mm}
\noindent
$\! \!${\foo dsuarez@dm.uba.ar}

\end{document}